\documentclass[oneside,12pt]{amsart}
\usepackage{
	amsmath,
	amssymb,
	amsthm,
	geometry,
	paralist,
	thmtools,
	tikz,
	tikz-cd,
	verbatim,
	relsize, 
	bbm, 
	dsfont,
	mathtools,
	enumerate
}
\usepackage[all]{xy}
\usepackage[T1]{fontenc}

\usepackage{caption}
\usepackage{subcaption}
\usetikzlibrary{angles}

\renewcommand\labelenumi{(\roman{enumi})}
\renewcommand\theenumi\labelenumi

\usepackage{adjustbox}

\tikzcdset{
	arrows = crossing over
}

\usepackage[maxbibnames=99,
maxalphanames=10, 
style=alphabetic]{biblatex}
\renewbibmacro{in:}{%
	\ifentrytype{article}{}%
	{\printtext{\bibstring{in}\intitlepunct}}}
\renewbibmacro*{url+urldate}{}
\renewbibmacro*{eprint}{}

\usepackage[colorlinks=true,allcolors=blue!80!black]{hyperref}

\usepackage{cleveref} 

\usepackage{dynkin-diagrams}

\usepackage[boxsize=.5em]{ytableau}

\usetikzlibrary{arrows.meta} 

\usepackage[textsize=scriptsize]{todonotes}
\presetkeys{todonotes}{inline}{}

\makeatletter
\providecommand\@dotsep{5}
\makeatother

\newcommand{\bbC}{\mathbb C}
\newcommand{\bbH}{\mathbb H}

\newcommand{\bbN}{\mathbb N}

\newcommand{\bbR}{\mathbb R}

\newcommand{\bbZ}{\mathbb Z}

\newcommand{\cB}{\mathcal B}
\newcommand{\cC}{\mathcal C}
\newcommand{\cD}{\mathcal D}

\newcommand{\cH}{\mathcal H}
\newcommand{\cI}{\mathcal I}

\newcommand{\cS}{\mathcal S}

\DeclareMathOperator{\Z}{\mathbb Z}

\DeclareMathOperator{\GL}{GL}

\DeclareMathOperator{\ra}{\rightarrow}

\DeclareMathOperator{\Sal}{Sal}
\DeclareMathOperator{\Tor}{Tor}

\declaretheorem[numberwithin=section]{theorem}
\declaretheorem[sibling=theorem]{lemma}

\declaretheorem[sibling=theorem]{conjecture}

\declaretheorem[sibling=theorem, style=remark]{remark}
\declaretheorem[sibling=theorem, style=definition]{definition}

\declaretheorem[sibling=theorem, style=definition]{exercise}

\newtheorem*{theorem*}{Theorem}

\crefname{lemma}{Lemma}{Lemma}
\crefname{corollary}{Corollary}{Corollary}
\crefname{theorem}{Theorem}{Theorem}
\crefname{definition}{Definition}{Definition}
\crefname{proposition}{Proposition}{Proposition}
\crefname{section}{Section}{Section} 
\crefname{construction}{Construction}{Construction}
\crefname{generalization}{Generalization}{Generalization}
\crefname{construction}{Construction}{Construction}
\crefname{notation}{Notation}{Notation}
\crefname{example}{Example}{Example}
\crefname{remark}{Remark}{Remark}
\crefname{fact}{Fact}{Fact}
\crefname{conjecture}{Conjecture}{Conjecture}
\crefname{motivation}{Motivation}{Motivation}  
\crefname{figure}{Figure}{Figure}  
\crefname{assumption}{Assumption}{Assumption}
\crefname{exercise}{Exercise}{Exercise}

\renewcommand{\comment}[1]{}

\newcommand{\G}{\Gamma}

\begin{document}
	
	\title[]{An introduction to the geometric and combinatorial group theory of Artin groups}
	
	\author[]{Rachael Boyd}
	\address{School of Mathematics and Statistics, University of Glasgow, Glasgow G12 8QQ, UK}
	\email{rachael.boyd@glasgow.ac.uk}
	\urladdr{https://www.maths.gla.ac.uk/~rboyd/} 
	
	\begin{abstract}
		We give a brief introduction to the geometric and combinatorial group theory of Artin groups. In particular we introduce the $K(\pi,1)$ conjecture for Artin groups and survey known results as of January 2024. These notes were written as companion notes for the MFO mini-workshop 2405a \emph{Artin groups meet triangulated categories} alongside Edmund Heng's notes \emph{Introduction to stability conditions and its relation to the} $K(\pi,1)$ \emph{conjecture for Artin groups}.
	\end{abstract}
	
	\maketitle
	
	
	\section{Introduction}
	
	There are many excellent survey articles about the geometric group theory of Artin groups. We start these notes with a non-exhaustive reference list for the reader:
	
	\begin{enumerate}
		\item[{\cite{CharneyProblems}}] Charney, \emph{Problems related to Artin groups}.
		\item[{\cite{GodelleParis2012}}] Godelle and Paris, \emph{Basic questions on Artin-Tits groups}.
		\item[{\cite{Paris2014}}] Paris, \emph{$K(\pi,1)$ conjecture for Artin groups}
		\item[{\cite{McCammondSurvey2017}}] McCammond \emph{The mysterious geometry of Artin groups}
	\end{enumerate}
	
	In these notes we will recap some of the salient definitions, results, methods, and open questions from the existing literature. We claim no originality in the mathematical content, aside from Conjecture~\ref{tor conjecture}. We include exercises and aim these notes at the non-expert.
	
	Let~$\Gamma$ be a finite labelled simplicial graph with vertex set~$S$, edge set~$E$, and edge labels~$m_{st} \geq 3$ or $m_{st}=\infty$, for $\{s,t\}\in E$. If there is no edge between vertices $s$ and $t$, let $m_{st}=2$. Then the Artin group~$A_{\Gamma}$ is given by the following presentation, where the relations are called \emph{generalised braid relations}.
	\[
	A_\Gamma =\langle S \mid \underbrace{sts\ldots}_{\text{length }m_{st}}=\underbrace{tst\ldots}_{\text{length }m_{st}} \,\, \forall \, s\neq t \in S \text{ such that }m_{st}\neq \infty \rangle.
	\]
	Note that when $m_{st}=2$ (there is no edge between $s$ and $t$ in $\Gamma$) then it follows that $s$ and $t$ commute. The Coxeter group~$W_\Gamma$ is given by adding the relations~$s^2=e$ for all $s \in S$. 
	
	There is a natural map $A_\G \to W_\G$ given by imposing the extra relations $s^2=e$, and the kernel of this map is the \emph{pure Artin group} $PA_\G$. The three groups therefore sit in the short exact sequence:
	
	\[
	PA_\G \hookrightarrow A_\G \to W_\G.
	\]
	
	There is also a canonical section to the projection map, $W_\G \to A_\G$, given by expressing a group element in the Coxeter group in terms of a minimal length word and viewing this as a word in the Artin group. Note that minimal words differ only by relations in $A_\G$, so the section does not depend on the choice of minimal word. We remark that there are in fact different choices one can make for this section, corresponding to the notion of \emph{linear braids}, but only the canonical section will appear in these notes.
	
	Finite Coxeter groups were classified by Coxeter in 1935 \cite{Coxeter1935}, motivated by the study of symmetry groups of regular polytopes. An Artin group $A_\G$ is called \emph{finite type} if the corresponding Coxeter group $W_\G$ is finite. We recall Coxeter's classification here, as it will provide many examples of groups for our exercises.
	
	\begin{theorem}[{Classification of finite Coxeter groups, Coxeter \cite{Coxeter1935}}]\label{prop:classification of finite Coxeter}
		A Coxeter group $W_\G$ is finite if and only if $\G$ is a disjoint union of finitely many of the following connected graphs.
		
		
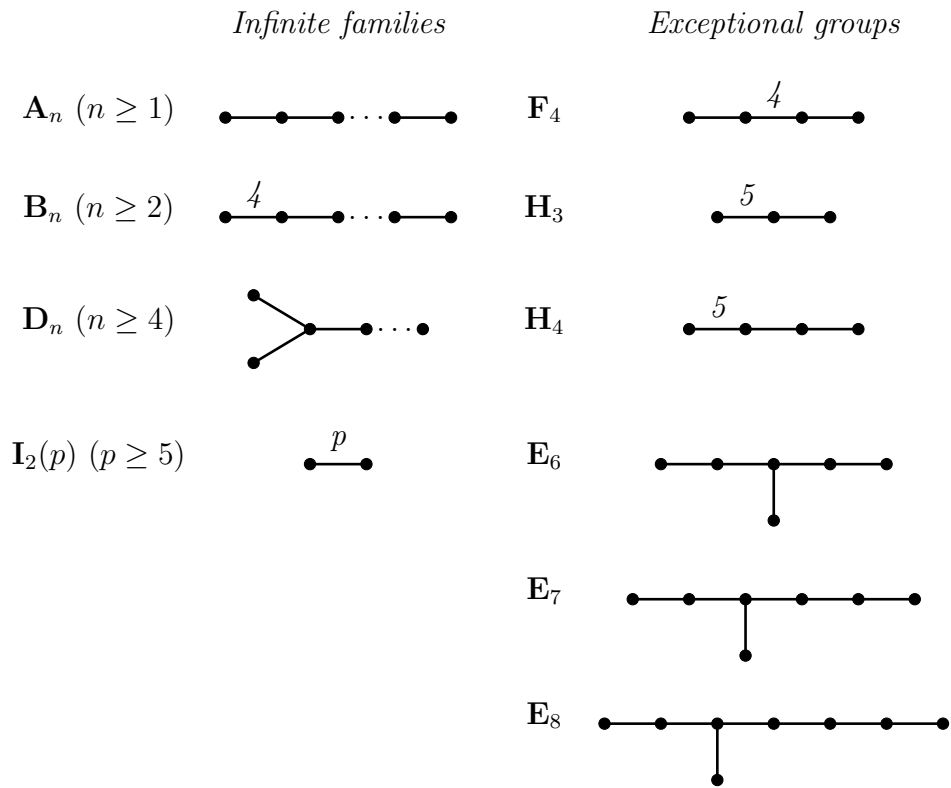
\centerline{	\xymatrix@R=3mm@C=2mm {	
				& \textrm{Infinite families} & & & \textrm{Exceptional groups}\\	
				{\bf A}_n \,\, (n\geq 1) & \begin{tikzpicture}[scale=0.15, baseline=0]
					\draw[fill= black] (5,0) circle (0.5);
					\draw[line width=1] (5,0) -- (15,0);
					\draw[fill= black] (10,0) circle (0.5);
					\draw[fill= black] (15,0) circle (0.5);		
					\draw (17.5,0) node {$\ldots$};
					\draw[fill= black] (20,0) circle (0.5);
					\draw[line width=1] (20,0) -- (25,0);				
					\draw[fill= black] (25,0) circle (0.5);
				\end{tikzpicture} 
				& & {\bf F}_4 & \begin{tikzpicture}[scale=0.15, baseline=0]
					\draw[fill= black] (5,0) circle (0.5);
					\draw[line width=1] (5,0) -- (20,0);
					\draw[fill= black] (10,0) circle (0.5);
					\draw[fill= black] (15,0) circle (0.5);
					\draw[fill= black] (20,0) circle (0.5);
					\draw (12.5,2) node {4};
				\end{tikzpicture}  \\
				{\bf B}_n \,\,(n\geq2) & \begin{tikzpicture}[scale=0.15, baseline=0]
					\draw[fill= black] (5,0) circle (0.5);
					\draw[line width=1] (5,0) -- (15,0);
					\draw[fill= black] (10,0) circle (0.5);
					\draw[fill= black] (15,0) circle (0.5);		
					\draw (17.5,0) node {$\ldots$};
					\draw[fill= black] (20,0) circle (0.5);
					\draw[line width=1] (20,0) -- (25,0);				
					\draw[fill= black] (25,0) circle (0.5);
					\draw (7.5,2) node {4};
				\end{tikzpicture} 
				& & {\bf H}_3 & \begin{tikzpicture}[scale=0.15, baseline=0]
					\draw[fill= black] (5,0) circle (0.5);
					\draw[line width=1] (5,0) -- (15,0);
					\draw[fill= black] (10,0) circle (0.5);
					\draw[fill= black] (15,0) circle (0.5);
					\draw (7.5,2) node {5};
				\end{tikzpicture}  \\
				{\bf D}_n \,\, (n\geq 4) & \begin{tikzpicture}[scale=0.15, baseline=0]
					\draw[fill= black] (5,-3) circle (0.5);
					\draw[line width=1] (5,-3) -- (10,0);
					\draw[line width=1] (5,3) -- (10,0);
					\draw[line width=1] (15,0) -- (10,0);
					\draw[fill= black] (5,3) circle (0.5);
					\draw[fill= black] (10,0) circle (0.5);
					\draw[fill= black] (15,0) circle (0.5);
					\draw (17.5,0) node {$\ldots$};				
					\draw[fill= black] (20,0) circle (0.5);
				\end{tikzpicture}  
				& & {\bf H}_4 & \begin{tikzpicture}[scale=0.15, baseline=0]
					\draw[fill= black] (5,0) circle (0.5);
					\draw[line width=1] (5,0) -- (20,0);
					\draw[fill= black] (10,0) circle (0.5);
					\draw[fill= black] (15,0) circle (0.5);
					\draw[fill= black] (20,0) circle (0.5);
					\draw (7.5,2) node {5};
				\end{tikzpicture}  \\
				{\bf I}_2(p) \,\, (p\geq 5)& \begin{tikzpicture}[scale=0.15, baseline=0]
					\draw[fill= black] (5,0) circle (0.5);
					\draw (7.5,2) node {p};
					\draw[line width=1] (5,0) -- (10,0);			
					\draw[fill= black] (10,0) circle (0.5);
				\end{tikzpicture} 
				& &{\bf E}_6 & \begin{tikzpicture}[scale=0.15, baseline=0]
					\draw[fill= black] (5,0) circle (0.5);
					\draw[line width=1] (5,0) -- (25,0);
					\draw[fill= black] (10,0) circle (0.5);
					\draw[fill= black] (15,0) circle (0.5);
					\draw[fill= black] (20,0) circle (0.5);		
					\draw[fill= black] (25,0) circle (0.5);			
					\draw[fill= black] (15,-5) circle (0.5);
					\draw[line width=1] (15,0) -- (15,-5);
				\end{tikzpicture}  \\ 
				& & & {\bf E}_7 & \begin{tikzpicture}[scale=0.15, baseline=0]
					\draw[fill= black] (5,0) circle (0.5);
					\draw[line width=1] (5,0) -- (30,0);
					\draw[fill= black] (10,0) circle (0.5);
					\draw[fill= black] (15,0) circle (0.5);
					\draw[fill= black] (20,0) circle (0.5);		
					\draw[fill= black] (25,0) circle (0.5);		
					\draw[fill= black] (30,0) circle (0.5);		
					\draw[fill= black] (15,-5) circle (0.5);
					\draw[line width=1] (15,0) -- (15,-5);
				\end{tikzpicture}   \\
				& & & {\bf E}_8 & \begin{tikzpicture}[scale=0.15, baseline=0]
					\draw[fill= black] (5,0) circle (0.5);
					\draw[line width=1] (5,0) -- (35,0);
					\draw[fill= black] (10,0) circle (0.5);
					\draw[fill= black] (15,0) circle (0.5);
					\draw[fill= black] (20,0) circle (0.5);		
					\draw[fill= black] (25,0) circle (0.5);		
					\draw[fill= black] (30,0) circle (0.5);		
					\draw[fill= black] (35,0) circle (0.5);		
					\draw[fill= black] (15,-5) circle (0.5);
					\draw[line width=1] (15,0) -- (15,-5);
				\end{tikzpicture}  \\ 
		}}
	\end{theorem}

	The generic definition of a Coxeter group was introduced by Tits in 1961 \cite{Tits1961} to accommodate infinite reflection groups associated to collections of hyperplanes: any Coxeter group can be realised as a discrete group generated by reflections on a finite dimensional vector space with respect to some inner product (the Coxeter group is finite precisely when this inner product is positive definite). We will define the hyperplane complements later. Artin groups were introduced by Brieskorn as the fundamental groups of complexified hyperplane arrangements, in relation to questions in algebraic geometry~\cite{Brieskorn1971}.

	Artin groups and Coxeter groups have been well studied due to their rich algebraic and geometric properties and connections to many fields of mathematics, such as representation theory, low dimensional topology, and monodromy theory. In geometric group theory, the theory of cube complexes lends itself to the study of \emph{Right Angled Artin groups} (RAAGs)---groups with only commuting ($m_{st}=2$) braiding relations and free relations ($m_{st}=\infty$)---see \cite{Wise2012} for a survey. In these notes we focus on the non-right-angled case, as the methods used to study RAAGs are somewhat distinct from general methods used to study Artin groups.

	\subsection{Open questions}
	Alongside the~$K(\pi,1)$ conjecture, which is the main topic of these notes, complete answers for some algebraic questions about Artin groups have remained out of reach with the current tool-kit, although many of them have been rephrased. The fundamental conjectures are that for general $\G$:
	\begin{enumerate}
		\item $A_\Gamma$ is torsion free.
		\item If~$A_\G$ cannot be written as a product of Artin groups, $A_\G$ has center $\bbZ$ when $A_\G$ is finite type, and trivial center otherwise.
		\item $A_\G$ has a solvable word (and conjugacy) problem.
	\end{enumerate}
	All of these conjectures are known to be true for certain families of Artin groups (e.g.~finite type Artin groups) but in general are unanswered, despite much work and interest in their direction. An overview can be found in \cite{GodelleParis2012} (where they prove that it suffices to prove the conjectures for $\Gamma$ ``free of infinity'') or \cite{McCammondSurvey2017}.
	
	We note here that a positive answer to the $K(\pi,1)$ conjecture implies that Artin groups are torsion free (we see later that the $K(\pi,1)$ conjecture being true exhibits a finite dimensional model for $BA_\G$ \cite{Salvetti1994}) and also answers the center conjecture \cite{JankiewiczSchreve2023}.
	
	There has been recent work into the \emph{isomorphism problem} for Artin groups (see, e.g.~\cite{vaskou2023, MartinVaskou2023}), which asks when two Artin groups are isomorphic (the corresponding question for Coxeter groups is still open). It is also an interesting open question to classify automorphisms of Artin groups, which has been done in very few cases e.g.~for RAAGs \cite{Droms1987}, and 2-dimensional Artin groups (defined later) \cite{Crisp2005}. Vaskou \cite{vaskou2023automorphisms} classifies automorphisms of large-type free-of-infinity Artin groups, and Bregman--Charney--Vogtmann have built an \emph{Outer space} for RAAGs \cite{BregmanCharneyVogtmann2023}.
	
	\subsection*{Acknowledgements}
	These notes were written for the Oberwolfach mini-workshop 2405a: \emph{Artin groups meet triangulated categories}.
	A companion set of notes can be found in \cite{Heng2024}.
	The corresponding Oberwolfach Report is \cite{MFOreport}.
	I would like to thank all participants of the MFO mini-workshop 2405a \emph{Artin groups meet triangulated categories} for their feedback on these notes. Particular thanks go to Andrea Bianchi, Edmund Heng, Viktoriya Ozornova, and Ailsa Keating for comments on a first draft.
	I would also like to thank my workshop co-organisers, Edmund Heng and Viktoriya Ozornova.
	Special thanks to Jon McCammond and Tony Licata for preparing lecture series for the mini-workshop.
	The workshop organisers would like to thank the Mathematisches Forschungsinstitut Oberwolfach (MFO) for the opportunity to organise a mini-workshop on this topic.
	The MFO and the workshop organisers would like to thank the National Science Foundation for supporting the participation of junior researchers in the workshop by the grant DMS-2230648, ``US Junior Oberwolfach Fellows''.
	
	\section{The $K(\pi,1)$ conjecture}
	In this section we state the $K(\pi,1)$ conjecture and introduce some various alternative forms.
	
	\subsection{Statement of the conjecture}

	\subsubsection{The canonical representation}
	
	For~$\G$ with vertex set~$S$, the Coxeter group $W_\G$ admits a faithful representation  $\rho\colon W_\G \to \GL(V)$ where $V$ is a real vector space of dimension~$|S|$. We describe this representation. 
	
	Let~$S=\{s_1,\ldots , s_n\}$ and let~$m_{ij}=m_{s_is_j}$. Let $V \cong \bbR^n$ with basis~$\{e_1, \ldots , e_n\}$.
	
	\begin{definition}
		Define a symmetric bilinear form $\cB$ on $V$ by
		\[
		\cB(e_i,e_j)=\begin{cases} 1 & \text{ if }i=j\\-\cos (\pi/m_{ij}) & \text{ if }m_{ij}\text{ finite }\\ -1 & \text{ if }m_{ij}\text{ infinite.}\end{cases}
		\]

		Using this we construct a map $\rho\colon W_\G \to \GL(V)$ by:
		\[
		\rho(s_i)= \sigma_i \qquad \sigma_i\colon V \to V; \sigma_i(v)=v-2\cB(e_i,v)e_i.
		\]
	\end{definition}

	\begin{exercise}
		Compute $\cB$ for the rank two Artin group ${\bf I}_2(p)$, and the affine Artin group $\widetilde{\bf A}_2$, which corresponds to $\G$ being a triangle with all edges labelled 3.
	\end{exercise}
	
	\begin{remark}
		The map $\sigma_i$ is linear, of order 2 and preserves the bilinear form $\cB$. Moreover $\sigma_i(e_i)=-e_i$ and the fixed points of $\sigma_i$ are precisely given by the hyperplane defined by the equation $\cB(e_i,v)=0$.
	\end{remark}
	
	\begin{exercise}
		Show that $\sigma_i\sigma_j$ has order $m_{ij}=p$ in the case of the rank two Artin group ${\bf I}_2(p)$.
	\end{exercise}

	In fact one can show that in general $\sigma_i\sigma_j$ has order $m_{ij}$, and the map $\rho$ extends to a homomorphism $\rho:W_\G \to \GL(V)$.
	
	\begin{theorem}[{\cite[Theorem 2.4]{Tits1961}}]
		The representation $\rho$ is faithful.
	\end{theorem}
	
	The group $W_\G$ can therefore be viewed as a subgroup of $\GL(V )$.

	\begin{remark}
		This symmetric bilinear form defines a matrix $B$ with $B_{ij}=\cB(e_i, e_j)$, sometimes called the \emph{Coxeter matrix} \cite{McCammondSurvey2017}. Since $\cB$ is real and symmetric, all of its eigenvalues are real. Then the spherical (finite) Coxeter groups are those for which all eigenvalues are positive ($\cB$ is positive definite), and the Euclidean Coxeter groups are those for which all eigenvalues are $\geq 0$ ($\cB$ is positive semi-definite). When the Coxeter matrix has at least one zero eigenvalue, the corresponding symmetric bilinear form is singular.
	\end{remark}
	
	\begin{definition}[Non-singular Tits cone]
		The faithful representation $\rho\colon W_\G \ra \GL(V)$ is generated by reflections and $W_\G$ acts properly discontinuously on a non-empty open cone $\cI \subseteq V$, called the \emph{Tits cone} with respect to $\rho$. When $B$ is non-singular the hyperplanes $H_i$ bound a closed simplicial cone $\cC$, and the Tits cone is defined to be the orbit $$
		\cI :=\bigcup_{w\in W_\G}w\cC.$$
	\end{definition}

	In the singular case, the hyperplanes $H_i$ no longer bound a simplicial cone, and Tits solved this problem by passing to the dual representation $\rho^*\colon W_\G \to \GL(V^*)$, where for $w\in W_\G, \phi \in V^*,$ and $v\in V$ we have 
	\[
	(\rho^*(w))(\phi)(v)=\phi(\rho(w^{-1})(v)).
	\]
	An alternative description is as follows. 
	With $\phi_i \in V^*$ denoting the standard dual basis such that $\phi_i(e_j) = \delta_{i,j}$, we have that $\rho^*(w)= (\rho(w)^{-1})^T \in \GL(V^*)$.
	This representation is faithful since it is dual to a faithful representation (and in fact it is easier to prove $\rho$ is faithful via proving that $\rho^*$ is faithful).
	
	In the dual representation $\rho^*\colon W_\G \ra \GL(V^*)$, the generators are once again acting by reflections, each of which fixes a hyperplane 
	\[
	H_i^*:=\{\phi \in V^* \, | \, \phi(e_i)=0\}
	\]
	and in this dual set up, the hyperplanes always bound a simplicial cone $$\cC^* := \{ \phi \in V^* \mid \phi(e_i) \geq 0\}.$$
	
	\begin{definition}[Tits cone; general case]
		Regardless of whether the matrix $B$ is singular or non-singular, the hyperplanes $H^*_i$ bound a closed simplicial cone $\cC^*$, and the \emph{Tits cone} with respect to $\rho^*$ is defined to be the orbit $$
		\cI :=\bigcup_{w\in W_\G}w\cC^*.$$
	\end{definition} 
	When $B$ is non-singular the two representations $\rho$ and $\rho^*$ are equivalent and so this definition is equivalent to that of the non-singular Tits cone.
	In fact, Charney and Davis further showed that the choice of $W_\Gamma$ as a linear reflection group is irrelevant when considering the homotopy type of the hyperplane complement (to be defined below) \cite[Section 2]{CharneyDavis1995}.
	
	Consider the set of all reflections in $W_\Gamma$, given by conjugates of the standard generators:
	\[
	R=\{wsw^{-1} \, | \, w\in W, s\in S\}.
	\]
	Then each reflection $r\in R$ satisfies that $\rho^*(r)\in \GL(V^*)$ fixes a hyperplane $H^*_r\subset \cI\subseteq V^*$. 
	\begin{remark}
		Note that the set of reflections in $W_\Gamma$ is \emph{not} the same as the set of elements of order two.
		Consider $\Gamma = {\bf A}_1 \amalg {\bf A}_1$ so that $W_\Gamma \cong \Z/2\Z \times \Z/2\Z$. This group is abelian, and so the reflections are simply the standard generators themselves. In particular, the element $(1,1)$ (which also has order 2) is not a reflection.
	\end{remark}
	\begin{exercise} 
		Show that under the representation $\rho$ (or $\rho^*$) defined previously, the element $\rho((1,1)) \in \GL(V)$ is a rotation by $\pi$ (consequently it has order 2).
	\end{exercise}
	Following \cite{GodelleParis2012}, we define the complexified hyperplane complement to be
	\[ \cH_{\G} = \cI \times V^* \setminus \cup_{r\in R} (H^*_r \times H^*_r).
	\]
	
	We call this the complexified hyperplane complement because in the non-singular case where the (non-singular) Tits cone satisfies $\cI\cong V$ (this happens when $W_\Gamma$ is finite), we have that
	\[
	\cH_\G \simeq V \otimes \bbC \setminus  \cup_{r\in R} (H_r \otimes \bbC).
	\]
	
	Note that in \cite{Paris2014}, $\cH_{\G}$ is defined to be $ \cI \times \cI \setminus \cup_{r\in R} (H^*_r \times H^*_r)$, which is homotopy equivalent to our definition by \cite[Lemma 2.1.1]{CharneyDavis1995}.
	
	This hyperplane complement is a connected manifold of (real) dimension $2|S|$ on which $W_\G$ acts freely and properly
	discontinuously. Van der Lek proved the following.
	
	\begin{theorem}[Van der Lek \cite{VanderLek1983}]
		The fundamental group of $\mathcal{H}_\G/W_\G$ is isomorphic to~$A_\G$.
	\end{theorem}

	Recall that for a group~$G$, a space $X$ is a $K(G,1)$ space if $X$ is a connected CW-complex, $\pi_1(X)\cong G$, and $\pi_i(X)=0$ for all $i \geq 2$.
	In this case, $X$ is also called a classifying space for $G$, denoted $BG$. Equivalenxtly, a CW-complex is a $K(G,1)$ space if and only if its universal cover is contractible.
	
	One of the forefront conjectures for Artin groups is the following, known as the \emph{$K(\pi,1)$ conjecture}.
	\begin{conjecture}[Arnol'd, Brieskorn, Pham, Thom]
		$\mathcal{H}_\G/W_\G$ is a classifying space for~$A_\G$, or $K(A_\G,1)$ space.
	\end{conjecture}

	\subsection{Alternative restatements of the conjecture}
	This conjecture has been rephrased in many ways, which we will now survey.

	The subgroup $A_T$ generated by a subset $T \subseteq S$ is called a special (or standard parabolic) subgroup of $A_\G$.  By a theorem of van der Lek \cite{VanderLek1983}, $A_T$ is isomorphic to the Artin group associated to the (full) subgraph of $\G$ spanned by $T$.

	\subsubsection{The Salvetti complex} 
	Charney and Davis \cite{CharneyDavis1995b}, and independently Salvetti \cite{Salvetti1994}, constructed a finite dimensional CW-complex, called the Salvetti complex and denoted by~$\operatorname{Sal}_\G$. They showed it is homotopy equivalent to~$\mathcal{H}_\G$, so  proving that this complex is aspherical also proves the~$K(\pi,1)$ conjecture.
	
	To build the Salvetti complex, we follow \cite{Paris2014} and \cite{Paolini2017}. We need some preliminaries on Coxeter group cosets.
	
	\begin{definition}
		Let $T\subseteq S$ and consider the special subgroup $W_T\subset W_\G$. We say that $w\in W$ is $T$-minimal if it is the unique minimal length representative in the coset $wW_T$. Such minimal length representatives always exist for Coxeter groups, see e.g.~\cite{Davis2008}.
	\end{definition} 
	
	Denote by $\cS^f$ the set $\{T\subseteq S \, | \, W_T \text{ is finite}\}$.
	
	\begin{definition}[Salvetti complex]
		The Salvetti complex is the geometric realisation of the poset given by endowing the set $W\times \cS^f$ with the following  partial order:
		\[
		(u, T) \leq (v, R)\text{ if }T \subseteq R, v^{-1}u \in W_R\text{ and }v^{-1}u\text{ is $T$-minimal}. 
		\]
	\end{definition}
	
	Charney and Davis, and independently Salvetti, proved the following. 
	\begin{theorem}[{\cite{Salvetti1994,CharneyDavis1995b}}]
		$\Sal_\G\simeq \cH_\G$.
	\end{theorem}
	
	The Salvetti complex $\Sal_\G$ comes with a natural action of the Coxeter group $W_\G$ by left multiplication and the equivalence $\Sal_\G\simeq \cH_\G$ is equivariant with respect to this action, so $\Sal_\G/W_\G \simeq \cH_\G/W_\G$ ({\cite{Salvetti1994,CharneyDavis1995b}}). The complex $\Sal_\G/W_\G$ admits a cell decomposition with one 0-cell, a one cell for each $s\in S$, and a $k$-cell for each $T\in \cS^f$ such that $|T|=k$. See \cite{Paris2014} and \cite{Paolini2017} for a precise account of this cell decomposition and some helpful figures.

	As an aside, when $T\subset S$, there is a natural embedding $\Sal_T\hookrightarrow \Sal_\G$ and Godelle and Paris prove the following interesting result:
	
	\begin{theorem}[{\cite[Theorem 2.2]{GodelleParis2012b}}]
		The natural embedding $\Sal_T\hookrightarrow \Sal_\G$ admits a retraction $\Sal_\G\rightarrow \Sal_T$.
	\end{theorem}
	
	We also provide an equivalent way of building the quotient of the Salvetti complex $\Sal_\G/W_\G$ using an Artin group action.
	
	Recall that there is a canonical section $W_\G \to A_\G$ (this is \emph{not} a group homomorphism). For $W_T$ a finite  special subgroup of $W_\G$, we denote by $\widehat{W}_T$ the image of $W_T$ under the section. Note that this is a subset of $A_\G$, not a subgroup.
	
	\begin{definition}[{\cite{CharneyProblems}}]
		Consider the geometric realisation $\Sal'_\G$ of the partially ordered set 
		\[
		\{a\widehat{W}_T \, | \, a \in A_\G, \, W_T \text{ finite}\}.
		\]
		Then the Artin group $A_\G$ acts freely and simplicially on $\Sal'_\G$. 
	\end{definition}
	
	\begin{lemma}
		$\Sal'_\G/A_\G \simeq \Sal_\G/W_\G.$
	\end{lemma}
	
	Note that in some later papers, e.g.~Salvetti and Paolini's proof of the $K(\pi,1)$ conjecture for affine Artin groups, the Salvetti complex is defined to be the quotient $\Sal_\G/W_\G$. In Paris's notes \cite{Paris2014} the quotient is denoted $\overline{\Sal_\G}$.

	\subsubsection{The Deligne complex}
	Building on the work of Deligne, Charney and Davis \cite{CharneyDavis1995} defined a cube complex whose vertices are given by cosets of finite type special subgroups of  $A_\G$.  Charney and Davis called this complex the modified Deligne complex, and it has since become known as the Deligne complex.  
	
	Recall that a cube complex is a space obtained by gluing Euclidean cubes of edge length 1 ($[0,1]^n$ for some $n\geq 0$) along subcubes ($[0,1]^k\subset [0,1]^n$ for some $n>k\geq 0$).
	
	\begin{definition}[The Deligne complex]\label{defn-Deligne complex}
		Let $A_\G$ be an Artin group. The \textit{Deligne complex}~$\cD_\G$ is the cube complex	with vertex set all cosets $gA_T$, such that $A_T$ is finite type ($T \in \cS^f$).  We partially order the vertices by inclusion of cosets. For any pair of vertices $gA_T \subset gA_{T'}$, the interval $[gA_T, gA_{T'}]$ spans a cube of dimension $|T' \smallsetminus T|$. 
	\end{definition}
	
	Note that the Deligne complex can also be constructed as a simplicial complex given by the geometric realisation of the poset $\{aA_T \, | \, a \in A_\G,  T\in \cS^f\}$. The cube complex structure is coarser than the simplicial one, and allows one to use cubical metrics such as the standard cubical metric induced from $[0,1]^n\subset \bbR^n$, and the Moussong metric.
	
	\begin{exercise}
		Draw the section of the Deligne complex corresponding to cosets with trivial representative, for the Artin groups of type~$\widetilde{\bf A}_2$, and ${\bf B}_3$. Give an example of a vertex with infinite valence in each of these complexes.
	\end{exercise}
	
	The action of $A_\G$ on its cosets by left multiplication induces a cocompact action by isometries of $A_\G$ on the Deligne complex~$\cD_\G$. Note, however, that this action is not proper as the stabiliser of a vertex $gA_T$ is the infinite subgroup $gA_Tg^{-1}$. 
	
	Charney and Davis proved the following:
	
	\begin{theorem}[{\cite{CharneyDavis1995}}]
		The Deligne complex $\cD_\G$ is homotopy equivalent to the universal cover of the Salvetti complex $\Sal_\G$.
	\end{theorem}
	
	We note than the analogously defined complex for Coxeter groups is the \emph{Davis complex}, but in this case the action of the Coxeter group is proper.

	\subsubsection{The extended Deligne complex}
	
	Godelle and Paris introduced a simplicial complex, that can also be given a cubical structure, which for the purpose of these notes we call the \emph{extended Deligne complex} \cite{GodelleParis2012b}. It extends the definition of Charney and Davis by replacing the family $\cS^f$ of finite type subsets of $S$ with a more general family $\cS$ of subsets. This family must satisfy the following properties:
	\begin{itemize}
		\item If $T \in \cS$ and $T'\subset T$, then $T' \in  \cS$,
		\item $A_T$ satisfies the $K(\pi, 1)$ conjecture for all $T\in \cS$,
		\item $\cS^f\subset \cS$.
	\end{itemize}
	A family $\cS$ of subsets of $S$ satisfying these conditions is called \emph{complete and $K(\pi,1)$}.
	
	\begin{definition}[The extended Deligne complex]\label{defn - extended Deligne complex}
		Let $A_\G$ be an Artin group, and $cS$ be a complete and $K(\pi,1)$ family of subsets of $S$. The \textit{extended Deligne complex}~$\cD_\G(\cS)$ is the geometric realisation of the poset of cosets $\{aA_T \, |\, a \in A_\G, \, T \in \cS\}$, ordered via inclusion.
	\end{definition}
	
	Godelle and Paris prove the following theorem.
	
	\begin{theorem}[{\cite[Theorem 3.1]{GodelleParis2012b}}]
		Let $\cS$ be a complete and $K(\pi, 1)$ family of subsets of $S$. Then $\cD_\G(\cS)$ is homotopy equivalent to the universal cover of the Salvetti complex $\Sal_\G$.
	\end{theorem}

	\subsubsection{The classifying space of the Artin monoid}
	
	The \emph{Artin monoid} $A_\G^+$ is the monoid given by the positive presentation
	\[
	A_\Gamma^+ =\langle S \mid \underbrace{sts\ldots}_{\text{length }m_{st}}=\underbrace{tst\ldots}_{\text{length }m_{st}} \,\, \forall s \neq t \in S \rangle^+.
	\]
	The group completion of $A_\Gamma^+$ is $A_\G$, and Paris showed that the inclusion map $A_\G^+ \hookrightarrow A_\G$ is injective \cite{Paris2002}.
	
	In 2006 Dobrinskaya proved that the quotient $\mathcal{H}_\G/W_\G$ has the same homotopy type as~$BA_\G^+$, the classifying space of the Artin monoid~\cite{Dobrinskaya2006}. 
	
	Recall the classifying space $BG$ of a discrete group $G$ can be built via the standard (or bar) resolution \cite[Section I.5]{BrownCohomology}, which has one vertex, and for $p\geq 0$, $p$-simplices correspond to tuples $[g_1| g_2 | \cdots | g_p]$ for $g_i \in G$. The first and final face map are given by forgetting~$g_1$ and $g_p$ and the intermediate face maps are given by forgetting a single bar and multiplying the group elements that it separated. The classifying space $BM$ of a monoid $M$ can likewise be built via the same bar construction, using monoid elements $m_i \in M$. The spaces obtained in this way are always connected and have the universal group associated with the monoid $M$ as the fundamental group. However, these spaces typically have non-vanishing higher homotopy groups. In fact, \cite{McDuff1979} shows that any connected space is weakly equivalent to $BM$ for some monoid $M$.
	
	The result by Dobrinskaya was later reproven by Ozornova \cite{Ozornova2017} and Paolini \cite{Paolini2017}, using discrete Morse theory.
	
	It follows that the~$K(\pi,1)$ conjecture is true for $A_\G$, if and only if the natural map~$BA_\G^+\to BA_\G$ is a homotopy equivalence. One can interpret this map concretely on bar constructions via the injection $A_\G^+\hookrightarrow A_\G$ inducing an injection on tuples.

	Interpreting work of Fiedorowicz \cite{Fiedorowicz1984} in the setting of Artin monoids, or alternatively considering Quillen’s Theorems A and B for the inclusion map $A^+_\G \hookrightarrow A_\G$, the $K(\pi,1)$ conjecture is equivalent to the following conjecture. 
	
	\begin{conjecture}[Boyd]\label{tor conjecture}
		$\Tor_*^{\bbZ A^+_\G}(\bbZ A_\G, \bbZ)=0$ for all $* \geq 1$.
	\end{conjecture}
	
	No current work on the conjecture explicitly states or utilises this rephrasing.

	\subsubsection{Overview}
	We give a diagrammatic depiction of some of the known~$K(\pi,1)$ conjecture equivalences in the diagram below. Here $\cS$ is any complete and $K(\pi,1)$ family of subsets of the standard generating set $S$. In particular when $\cS=\cS^f$, $\cD_\G(\cS)$ is the Deligne complex $\cD_\G$. Finding a proof which confirms any question mark shown in the diagram would in turn prove the~$K(\pi,1)$ conjecture.
	\[
	\xymatrix@C=7em{
		\cD_\G(\cS) \simeq\widetilde{\Sal_\G}\overset{\text{\large{?}}}{\simeq}\ast \qquad\ar[d]&&\\
		\operatorname{Sal}_\G/W_\G\ar[r]^{\simeq}&\mathcal{H}_\G/W_\G \ar[r]^{K(\pi,1) \text { conj.}}_{\simeq \text{\large{?}}}& BA_\G\\
		&BA^+_\G \ar[u]_{\simeq} \ar[ru]_{\simeq \text{\large{?}}}\\
	}
	\]

	\subsection{Progress on the conjecture}
	We now give an overview of the families of Artin groups for which the $K(\pi,1)$ conjecture has been solved.\footnote{This survey was written in January 2024 -- further developments since then are not covered in this survey.}
	
	\subsubsection{Finite type Artin groups}
	Recall that an Artin group is finite type, or spherical type, if the associated Coxeter group is finite. Deligne proved the~$K(\pi,1)$ conjecture for finite type Artin groups \cite{Deligne1972}, using a complex which would later be modified by Charney and Davis, and dubbed the Deligne complex (\cref{defn-Deligne complex}). Deligne's proof utilises the Garside structure on the Artin monoid that we introduce in the next section, and a `union of chambers' argument, which we will also give a general framework for. In his notes on the $K(\pi,1)$ conjecture Paris reproves Deligne's result, by showing that the universal cover of the Salvetti complex is contractible. This proof is very accessible and we recommend it -- it also follows a union of chambers argument and utilises properties of the monoid in the finite type case \cite{Paris2014}.
	
	\subsubsection{FC-type Artin groups}
	An Artin group is called \textit{FC-type} if any subset $T\subseteq S$ that does not contain $s,t$ with $m_{st}=\infty$ generates a finite type Artin group, i.e.~$T\in \cS^f$.   FC-type Artin groups were originally defined by Charney and Davis in \cite{CharneyDavis1995}. They are precisely those for which the standard cubical metric on the Deligne complex is CAT(0).``FC" stands for Flag Complex, which comes from Gromov's flag condition on cube complexes: if the links of all vertices are flag complexes, then the cube complex is CAT(0).  The Deligne complex satisfies this condition if and only if $A_\Gamma$ is FC-type, by definition. In addition to proving the $K(\pi,1)$ conjecture is true in this case, the Deligne complex has been used to show that FC-type Artin groups have solvable word problem, are torsion-free and have finite virtual cohomological dimension, among other properties \cite{Altobelli1998, CharneyDavis1995, Godelle2007}.  
	
	\subsubsection{2-Dimensional and locally reducible Artin groups}
	
	An Artin group is said to be \emph{2-dimensional} if all finite type $T\subseteq S$ satisfy that $|T|\leq 2$ (type ${\bf A}_1$ or ${\bf I}_2(p)$). In particular this means that the Salvetti complex and the Deligne complex are 2-dimensional. For these groups the $K(\pi,1)$ conjecture was proved by Charney and Davis \cite{CharneyDavis1995}, and Charney extended the class to those Artin groups for which the finite type special subgroups are either those from the 2-dimensional case, or the braid group on 4 strands (type ${\bf A}_3$), \cite{Charney2004}. This extended family are called \emph{locally reducible} Artin groups. Both of these results were proved by showing that the Moussong metric on the Deligne complex satisfies non-positive curvature properties (CAT(0) in the 2-dimensional case, and CAT(1) in the extended case).
	
	\subsubsection{Large type Artin groups}
	An Artin group $A_\G$ is said to be \emph{of large type} if every $m_{st}$ appearing in the presentation is $\geq 3$, i.e.~no two standard generators commute. Hendriks proved the $K(\pi,1)$ conjecture for Artin groups of large type \cite{Hendriks1985}. Note that large type Artin groups are 2-dimensional so this result is recovered by the previous case.
	
	\subsubsection{Euclidean Artin groups}
	The combinatorial structure of Artin groups was first studied for the braid groups by Garside \cite{Garside1969}, who found an elegant solution to the word problem. This was extended to finite type Artin groups by Brieskorn and Saito~\cite{BrieskornSaito1972} and the combinatorial consequences have played a major role in the study of finite type Artin groups: for instance this structure can be used to prove the~$K(\pi,1)$ conjecture. The notion of a {Garside structure} was later introduced and explored in detail in \cite{DehornoyParis1999} (see also \cite{DDGKM}) -- we will define this later in these notes.
	
	More recently, a {dual Garside structure} was defined for affine Artin groups \cite{McCammondSulway}, and this was used by Paolini and Salvetti \cite{PaoliniSalvetti2021} to prove the~$K(\pi,1)$ conjecture for the class of affine Artin groups. Paolini has an accessible set of notes on their proof \cite{paolini2021dual} which also outlines the dual approach to the $K(\pi,1)$ conjecture in general. With Delucchi they extended their methods to prove the conjecture for all Artin groups of rank 3 \cite{DelucchiPaoliniSalvetti2022}.
	
	\subsubsection{Recent progress}
	Haettel and Huang recently proved the $K(\pi_1)$ conjecture for a family of Artin groups by showing that~$A_\Gamma \times \mathbb{Z}$ exhibited a Garside structure, and using this structure to deduce properties of~$A_\Gamma$ \cite{HaettelHuang}. This family includes, and is built from, Artin groups of \emph{cyclic type}. They define an Artin group to be of cyclic type if its defining graph~$\G$ is a cycle, and any special subgroup is finite type. Their new examples of Artin groups satisfying the  $K(\pi,1)$ conjecture contains some $A_\G$ for which $W_\G$ acts on the hyperbolic spaces $\bbH_3$ or $\bbH_4$.
	
	At approximately the same time, Huang released a preprint \cite{huang2023} proving the $K(\pi,1)$ conjecture for $A_\G$ where $\G$ is a tree satisfying certain properties, or $\G$ has a cyclic subgraph and satisfies certain properties. For instance, he proves the conjecture when $\G$ is a tree and has a collection of edges $E$ such that the components of $\Gamma\setminus E$ are either spherical, affine or locally reducible (so in particular the conjecture is already known). Huang proves the conjecture by introducing and studying the homotopy type of \emph{relative Artin complexes} and \emph{folded Artin complexes}. We briefly discuss the relative Artin complex.
	
	\begin{definition}
		The \emph{Artin complex} $\Delta_\G$ is the simplicial complex with vertex set the set of cosets~$\{aA_T \, | \, a \in A_\G, \, T=S\setminus \{s\} \text{ for } s\in S \}$.  We say a vertex is of \emph{type} $\hat{s}$ if it is a coset of $A_{S\setminus s}$. A collection of vertices span a simplex if the associated cosets have nonempty common intersection. 
	\end{definition}
	
	Under the assumption that $S\setminus \{s\}$ satisfies the~$K(\pi,1)$ conjecture for all~$s\in S$, the Artin complex agrees with the extended Deligne complex $\cD_\G(\cS)$ from \cref{defn - extended Deligne complex} with~$\cS$ the family of proper subsets of $S$. Huang's proof strategy starts from the fact that showing the extended Deligne complex is contractible proves the conjecture. He introduces the relative Artin complex which we define here for the interested reader.

	\begin{definition}[{\cite[Definition 2.3]{huang2023}}]
		For $A_\G$ an Artin group, let $S'\subset S$. The \emph{relative Artin complex} $\Delta_{S, S'}$ is the induced complex of $\Delta_\G$ spanned by vertices of type~$\hat{s}$ for $s\in S'$.
	\end{definition}
	
	Huang shows that the following conjecture is equivalent to the~$K(\pi,1)$ conjecture \cite[Corollary 7.4]{huang2023}. An Artin group~$A_\G$ is called \emph{almost spherical} if $\G$ is free of infinity (no $m_{st}=\infty$), $A_\G$ is not finite type, but $A_{T}$ is finite type for all proper subsets $T\subset S$. 
	
	\begin{conjecture}[{\cite[Conjecture 2.4]{huang2023}}]
		Suppose $\G$ is a complete graph with vertex set $S$ and $T\subset S$ is such that $A_T$ is almost spherical. If $\G$ and the subgraph spanned by $T$ are both connected and free of infinity, then $\Delta_{S, S'}\simeq *$.
	\end{conjecture}
	
	Note that in both of these papers the graph $\Gamma$ we used to define $A_\G$ is denoted $\Lambda$ and $\Gamma$ is used for a different graph convention.

	\subsubsection{Conjecturette} 
	The $K(\pi,1)$ \emph{conjecturette} \cite{EliasWilliamson2017} states that~$\pi_2(\mathcal{H}_\Gamma/W_\Gamma)=0$.  Elias and Williamson originally claimed that the conjecturette was proved in~\cite[\S 6]{DigneMichel}, but they have since retracted this claim and hence the conjecturette remains open.
	
	\section{Some techniques used in the study of Artin groups}
	
	In this section we give an overview of some of the basic theory used to prove some of the results discussed above.
	
	\subsection{Garside theory}
	
	Recall that the combinatorial structure of Artin groups was first studied for the braid groups by Garside \cite{Garside1969} and extended to finite type Artin groups by Brieskorn and Saito~\cite{BrieskornSaito1972}. Since the combinatorial consequences have played a major role in the study of finite type Artin groups, Dehornoy and Paris defined a group with such a structure to be a \emph{Garside group}. This led to the development of \emph{Garside theory} which asked the question: What are features of groups (and more generally monoids or categories) which allow us to carry over (analogs of) Garside's proof? After several different but related notions in the literature, the one in \cite{DDGKM} became somewhat standard. We define the structure and list some consequences below.
	
	\begin{definition}[Garside monoid {\cite[Definition 2.1]{DehornoyParis1999}}]
		A \emph{Garside monoid} is a pair $(M, \Delta)$ where $M$ is a monoid and $\Delta \in M$ such that the following hold.
		\begin{enumerate}[(i)]
			\item M is left- and right-cancellative.
			\item There exists $\ell \colon M \to \bbN$ satisfying $\ell(ab) \geq \ell(a) + \ell(b)$, and $\ell(a)=0 \Rightarrow a=e$ (this induces a partial order on $M$).
			\item Any two elements of $M$ have a left- and a right-lcm (lowest common multiple) and a left- and a right-gcd (greatest common divisor).
			\item $\Delta$ is a \emph{Garside element} of M, i.e.~the set of left-divisors of~$\Delta$ is equal to the set of right-divisors, and this set is a generating set for $M$.
			\item The set of all divisors of $\Delta$ in $M$ is finite.
		\end{enumerate}
	\end{definition}
	
	Given a Garside monoid, we can construct exactly one \emph{Garside group}, which is the enveloping group (or group-completion) of the monoid. It is then a property of the monoid that this group-completion can in fact be constructed by left fractions, and that the monoid injects into the group-completion.
	
	\begin{definition}[Garside group {\cite[Definition 2.3]{DehornoyParis1999}}]
		$G$ is a \emph{Garside group} if it is the group completion of a Garside monoid $M$. In this case, it follows from $M$ being Garside that:
		\begin{itemize}
			\item $M \subseteq G$,
			\item for all $g \in G$, theres exists $a, b \in M$ such that $g=a^{-1}b$.
		\end{itemize}
		We say that $G$ is the \emph{group of left-fractions} of $M$.
	\end{definition}
	
	A monoid is a \emph{quasi-Garside monoid}, if it satisfies the properties of a Garside monoid except the final one: $\Delta$ is allowed to have an infinite number of divisors. In this case the group of fractions is called a \emph{quasi-Garside group}.

	Let $A_\G^+$ be an Artin monoid with generating set $S$ and let $R_\G$ be the set of relations determined by the defining graph $\G$. Consider the set $\{S\}^*$ of words written in the alphabet $S$. We often consider these words modulo the relations $R_\G$. This set is in one-to-one correspondence with the set of elements in the Artin monoid~$A_\G^+$. If $w$ denotes a word in $\{S\}^*$, let $\overline{w}$ denote its image in $\{S\}^*/R_\G \cong A_\G^+$. We can endow the set $\{S\}^*$ with a partial order $\leqslant_R$, where we say two words $v$ and $w$ in $\{S\}^*$ satisfy $v\leqslant_R w$ if there exists $z\in \{S\}^*$ such that $w=zv$ i.e.~\emph{$v$ appears on the right of $w$}. This partial order can be passed to the Artin monoid: we say two elements $\overline{v}$ and $\overline{w}$ in $\{S\}^*/R_\G$ satisfy $\overline{v}\leqslant_R \overline{w}$ if there exists $z\in \{S\}^*$ such that $\overline{w}=\overline{zv}$. The partial order $\leqslant_L$ is defined similarly, where $v\leqslant_L w$ if there exists $z\in \{S\}^*$ such that $w=vz$ i.e.~\emph{$v$ appears on the left of $w$}. Under the isomorphism $\{S\}^*/R_\G \cong A_\G^+$ we can write $a$ and $b$ in $A_\G^+$ satisfy $a\leqslant_R b$ if there exists $c\in A^+_\G$ such that $b=ca$. The set of $b\in A_\G^+$ satisfying $b\leqslant_R a$ are called the \emph{right divisors} of $a\in A_\G^+$, and similarly for the set of \emph{left-divisors.}
	
	\begin{exercise}
		Determine the exact definitions for the notions of least common multiple and greatest common divisor for $A_\G^+$ with the partial order $\leqslant_R$ or $\leqslant_L$.
	\end{exercise}

	It was shown by Brieskorn and Saito that finite type Artin monoids have a Garside structure, which we now introduce. Recall that a finite Coxeter group $W_\G$ has a unique longest element, and let $\Delta$ be the image of this element in $A_\G$ under the canonical section $W_\G \to A_\G$. This lies in~$A_\G^+$, since the entire image of $W_\G$ under the canonical section does. Then the pair $(A_\G^+, \Delta)$ is a Garside monoid since:
	\begin{enumerate}[(i)]
		\item \cite[Proposition 2.3]{BrieskornSaito1972} $A_\G^+$ is left and right cancellative, i.e.~if $ac=bc$ or $ca=cb$ in $A_\G^+$, then we can cancel the $c$ term on both sides and conclude that $a=b$.
		\item The \emph{length function} $\ell\colon A_\G^+ \to \bbN$ which maps an element~$a\in A^+_\G$ to the length of a word representing $a$ (in the standard generators) is well defined, since all relations in~$A_\G^+$ are the same length on either side. This function satisfies that $\ell(ab)=\ell(a)+\ell(b)$, and the only monoid element which evaluates to 0 is the identity element $e \in A^+_\G$.
		\item \cite[Proposition 4.2 and Theorems 5.5 and 5.6]{BrieskornSaito1972} In a finite type Artin group, every set of elements has a least common multiple and greatest common divisor.
		\item $\Delta$ is a Garside element, and its set of left-divisors equals its set of right-divisors. This set is exactly the image of the canonical section $W_\Gamma\to A_\Gamma^+$ and hence always contains the standard generating set $S$, and therefore generates $A_\G^+$. In particular, this means for every $s\in S$ there exists $a,b\in A_\G^+$ such that $\Delta=s a = b s$.
		\item Since $W_\G$ is finite, the set of divisors of $\Delta$ is therefore finite.
	\end{enumerate}
	
	In contrast to item (iii), in a generic Artin group greatest common divisors always exist and least common multiples exist only when a common multiple exists \cite[Proposition 4.1]{BrieskornSaito1972}. 
	
	It follows that in a finite type Artin group $A_\G$, all group elements $g$ can be written as $g=ab^{-1}$ for $a, b \in A_\G^+$. Moreover, for each $g$ there exists an integer $k$ and element $a\in A_\G^+$ such that $g=\Delta^ka$.
	
	Moreover, it is an interesting fact that given a finite type Artin monoid $A_\G^+$ with generating set $S$ and Garside element $\Delta$, there exists a permutation $\sigma\colon S \to S$ such that for all $s\in S$, $\Delta s=\sigma(s)\Delta$.
	
	\begin{exercise}
		For the monoids associated to the finite type Artin groups of type ${\bf I}_2(p)$, and ${\bf A}_3$, determine the Garside element $\Delta$, and permutation $\sigma\colon S \to S$.
	\end{exercise}
	
	Given an Artin monoid~$A_\G^+$, with generating set~$S$, we denote its Garside element (if it exists) by~$\Delta_S$. For $T\subseteq S$, we denote the Garside element for the submonoid $A_T^+$ by $\Delta_T$. The following lemma gives rise to a useful normal form for Artin monoids.
	
	\begin{lemma}
		If $a\in A_\G^+$ has $T\subseteq S$ as its set of right-divisors of length 1, then $A_T^+$ is a finite type Artin monoid with Garside element~$\Delta_T$, and there exists $b\in A_\G^+$ such that $a=b\Delta_T$.
	\end{lemma}
	
	We can use this fact first to write $a=b\Delta_{T_0}$, then to write $b=c\Delta_{T_1}$ where $T_1\subset S$ is the set of right-divisors of $b$ with length 1, and so on until we reach the unique normal form:
	\[
	a=\Delta_{T_k}\ldots \Delta_{T_1}\Delta_{T_0}.
	\]
	Here $a$ has been written in a normal form using the alphabet of Garside elements corresponding to finite type Artin monoids, but $A_\G^+$ does not have to be a finite type Artin monoid itself. This normal form has been key to many proofs for Artin monoids, but passing from monoids to groups is still very difficult without a (quasi-) Garside structure.

	Dehornoy has a series of papers \cite{Dehornoy1, Dehornoy2, Dehornoy3, Dehornoy4} outlining a program to export the essence of Garside theory to general Artin groups. As we have seen, for finite type Artin groups, which are Garside groups, one can write every element as~$g=a^{-1}b$ for~$a$ and~$b$ in the associated monoid. This is not true for general Artin groups, and Dehornoy and coauthors study rewriting systems on $A_\Gamma$ arising from writing every element as a \emph{multifraction} $g=a_1^{-1}b_1\ldots a_n^{-1}b_n$. Dehornoy reformulates the word problem for $A_\Gamma$ using multifractions, and formulates several conjectures. The most basic of these is the following: suppose an element $g\in A_\Gamma$ can be written as~$g=a^{-1}b$ with no cancellation, then are~$a$ and~$b$ unique? Techniques in the study of $\mathrm{Stab}(D_{\Gamma})$ (see the companion notes \cite{Heng2024}) have shown promise to solve this basic conjecture, and thus it is conceivable that these techniques in conjunction with Dehornoy's program may shed further light on the word problem for Artin groups.
	
	\subsection{Dual Garside structures}
	
	\begin{definition}
		A \emph{Coxeter element} for an Artin group~$A_\G$ with finite generating set $S=\{s_1, \ldots , s_n\}$ is an element represented by a word $\sigma(s_1)\sigma(s_2)\ldots \sigma(s_n)$ where $\sigma$ is a permutation of the $n$ elements of $S$. For example, $s_1\cdots s_n$ is a Coxeter element for $A_\G$.
	\end{definition} 
	
	Given an Artin group $A_\G$ and Coxeter element $\delta$ one can define a \emph{dual Artin group} (e.g.~\cite{paolini2021dual}). When $A_\G$ is finite type, this group is both isomorphic to the Artin group, and Garside, with Garside element given by the Coxeter element~$\delta$. Finite type Artin groups therefore have two (different) Garside structures, corresponding to the standard presentation and dual presentation \cite{Bessis03}.
	
	In general, it is conjectured that the dual Artin group is isomorphic to the Artin group. This conjecture has been shown to be true for rank 3 Artin groups \cite{DelucchiPaoliniSalvetti2022} and affine Artin groups \cite{Digne2006, Digne2012, McCammond2015, McCammondSulway}.
	
	The dual Garside structure for affine Artin groups introduced by McCammond and Sulway \cite{McCammondSulway} was used by Paolini and Salvetti \cite{PaoliniSalvetti2021} to prove the~$K(\pi,1)$ conjecture for the class of affine Artin groups.

	\subsection{Shellability and union of chambers}

	In \cite{Davis2008}, Davis used a union of chambers argument to prove that the Davis complex $\Sigma_W$ associated to a Coxeter group is contractible. He did this by showing that the Davis complex is an example of a \emph{basic construction}. For a second example, Hepworth's high connectivity results relating to homological stability for Coxeter groups \cite{HepworthCoxeter} also used such an argument. In \cite{Paris2014}, Paris used a union of chambers argument to show that the universal cover of an analogue of the Salvetti complex for certain Artin monoids is contractible, proving the $K(\pi,1)$ conjecture for finite type Artin groups. Loosely, the argument consists of decomposing the complex into a union of high dimensional \emph{chambers} and considering how connectivity changes as they are glued together to build up the complex.

	We give an overview of the argument in the case that we have a simplicial complex $C$ of dimension $n$ that we wish to prove is contractible, or $(n-1)$-connected. The argument works in the setting that $C$ is a union of $n$ dimensional simplices over some (possibly infinite) indexing set $I$:
	\[
	C=\bigcup_{i\in I} \sigma^n_i
	\]
	and to each simplex there is a canonical way to associate a natural number, i.e.~there exists a map $l\colon I\to \bbN$, such that only one simplex -- $\sigma_{i_0}^n$ -- has index which evaluates to $0$. This induces a partial ordering on the $\sigma^n_i$, inherited from $\leq$ on $\bbN$. For instance, if $C$ is the geometric realisation of a poset of group cosets, the indexing set may be the group, or a set of coset representatives, and a well defined length function may act as the map to $\bbN$. Alternatively, $C$ may be the orbit of an $n$-dimensional simplex under a group action, again with the group acting as the indexing set.
	We filter $C$ by the natural numbers as follows:
	
	\begin{definition}
		For $k$ in $\mathbb{N}$ we define $C(k)$ as follows:
		\begin{equation*}
			C(k)=\bigcup_{\substack{i \in I, \\ l(i)\leq k}} \sigma^n_i.
		\end{equation*}
	\end{definition}
	
	We remark that $C(0)=\sigma^n_{i_0}$, $C(k) \subseteq C(k+1)$, and the complex $C$ is given by $\lim_{k\to \infty}C(k)$. The union of chambers argument relies on proving the following two claims:
	
	\begin{enumerate}
		\item[Claim (A):]  If $l(i)=k+1$ then $ \sigma_i^n \cap C(k)$ is a non-empty union of top dimensional faces of $\sigma_i^n$.
		\item[Claim (B):]  If $l(i)=l(j)=k+1$ and~$i\neq j \in I$ then $\sigma_i^n \cap \sigma_j^n \subseteq C(k)$.
	\end{enumerate}
	
	These claims must be verified for the individual case being considered. Assuming these claims, we will now show how to prove, that $C$ is either contractible or $(n-1)$-connected.
	
	We inductively build $C$ by increasing $k$ in $C(k)$, and observe how the homotopy type changes. We start at $C(0)=\sigma^n_{i_0}$, which is contractible. At each step we build up from $C(k)$ to $C({k+1})$ by adding the set of simplices:
	$$
	\bigcup_{\substack{i\in I, \\ l(i)= k+1}} \sigma^n_i.
	$$
	Then Claim (A) implies that $C(k)\cup \sigma^n_i$ is a non-empty union of facets of $\sigma_i^n$. Therefore passing from $C(k)$ to  $C(k)\cup \sigma^n_i$, the effect on the homotopy type is either: 
	\begin{itemize}
		\item it doesn't change upon adding the simplex (if not all facets of $\sigma_i^n$ are in the intersection), or
		\item the homotopy changes and this change is described by the possible addition of an $n$ sphere (if all facets are in the intersection, i.e.~$ \sigma_i^n \cap C(k)=\partial\sigma_i^n$ ).
	\end{itemize} 
	Claim (B) now implies that adding the entirety of the above union of simplices to $C(k)$ concurrently only changes the homotopy type in the sense that the individual simplices change it, since each two simplices intersect within $C(k)$. Therefore at each stage $C(k)\to C(k+1)$ we change the homotopy type by at most the addition of $n$-spheres and it follows from passing to the limit that $C$ is either contractible, or $(n-1)$-connected.
	
	The function~$l\colon I \to \bbN$ gives a partial order on the top dimensional simplices of~$C$. By~Claim (B), any linear extension of this partial order to a total order will still satisfy~Claim (A). In this case, the ordering is called a \emph{shelling} (see~\cite{Bjorner1992}), which we know to be highly connected, giving an alternative proof to the previous discussion.
	
	This proof outline can be easily altered. For example, we can start with $C(0)$ being a $m$-connected complex for some $m$. We can also choose the fundamental building block to be more complicated than a single simplex. For example, it could be the fundamental domain for a group or monoid action on a cube complex or simplicial complex.
	
	\printbibliography
	
%
	
\end{document}